\documentclass{article}
\usepackage{epsfig}
\usepackage{latexsym}
\usepackage{amsmath}
\usepackage{amsthm}
\usepackage{amssymb}
\usepackage{graphicx}
\usepackage{mdwlist}

\usepackage{thesis}
\usepackage{hodge}
\usepackage{eufrak}
\begin{document}
\title{Evaluating tautological classes using only Hurwitz numbers}
\author{Aaron Bertram, Renzo Cavalieri and Gueorgui Todorov}

\maketitle
\begin{abstract}
Hurwitz numbers count ramified covers of a Riemann surface
with prescribed monodromy. As such, they are purely combinatorial objects.
Tautological classes, on the other hand, are distinguished classes in the
intersection ring of the moduli spaces of Riemann surfaces of a given
genus, and are thus ``geometric.'' Localization computations in
Gromov-Witten theory provide non-obvious relations between the two.
This paper makes one such computation, and shows how it leads to
a ``master'' relation (Theorem \ref{rel}) that reduces the ratios of certain
interesting tautological classes to the pure combinatorics of Hurwitz numbers. As a corollary, we obtain a purely combinatorial proof of a theorem of Bryan and Pandharipande, expressing in generating function form classical computations by Faber/Looijenga (Theorem \ref{bpa}).
\end{abstract}
\section*{Introduction}

Clever applications of the Atiyah-Bott localization theorem in the context of 
Gromov-Witten theory have generated volumes of enumerative data as well as
many insights
into the structure of the
tautological rings of the moduli spaces of curves (\cite{fp:hiagwt}, \cite{fp:rmatc}, \cite{gjv:gwpophnahh}, \cite{gjv:spolgc}). The general idea is to exploit
torus actions on a target manifold (often just $\mathbb P^1$) to obtain torus actions on 
the moduli spaces of   
stable maps to the target. An analysis of the fixed loci for the torus action  then produces intersection numbers and relations among tautological classes.

Recent work (\cite{gv:taut}, \cite{gv:hnavl},\cite{gv:rvl}, \cite{i:rittr}) has shown that the same idea, when 
applied to moduli of ``relative'' stable maps reveals even more information.
In this paper, we follow ideas of the second author in \cite{r:dd2}, applying localization 
in the context of the ultimate relative stable map spaces to $\mathbb P^1$, the admissible cover spaces 
of \cite{acv:ac}, to tie together enumerative data
of apparently quite different natures by means of one succinct formula.
\begin{theorem}
\label{rel}
For each fixed integer $d\geq1$,
\begin{eqnarray}
	\mathcal{CY}(u)=\frac{(-1)^{d-1}d!}{d^{d}} \hf{(d)}(u) \mbox{e}^{d\mathcal{D}(u)},
\end{eqnarray}
where
\begin{basedescript}{\desclabelstyle{\pushlabel}\desclabelwidth{1.1cm}}
\item[$\mathbf{\hf{(d)}(u)}$] is the generating function for degree $d$, one-part simple Hurwitz numbers (properly defined in section \ref{hurwitz}).
\item[$\mathbf{\mathcal{D}(u)}$] is the generating function for $\la$ Hodge integrals on the moduli spaces $\mathcal{A}^g_{dd}$ (a.k.a. the evaluation of $[\mathcal{A}^g_{dd}]$ in $R^{g-2}(\mathcal{M}_g))$(discussed in the next paragraph and in section \ref{la}).
\item[$\mathbf{\mathcal{CY}(u)}$] is the generating function for the fully ramified Calabi-Yau cap invariants of degree $d$ maps in \cite{r:tqft}, briefly described in section \ref{CYcap}.
\end{basedescript} 
\end{theorem}
This new formula combines with other known
formulas to explicitly reduce the 
enumerative invariants involved in the 
local Gromov-Witten theory of curves (\cite{bp:tlgwtoc}) to 
the ``pure combinatorics'' of single Hurwitz numbers. 

\medskip

One class of invariants appearing in Theorem \ref{rel} is worth singling out. The tautological 
ring $R^*({\mathcal{M}_g})$ of the moduli stack of curves of genus $g$ 
is the subring of the Chow ring (with rational coefficients) 
generated by Mumford's ``kappa'' classes (see \cite{m:taegotmsoc}). It is conjectured in \cite{f:acdottr}
that this ring satisfies Poincar\'e duality with a (one-dimensional) socle in degree $g-2$. We focus on a particular family of  degree $g-2$ tautological classes, which arise naturally from Hurwitz theory.
Consider the $(2g-1)$-dimensional ``Hurwitz spaces:'' 
$$        \mathcal{A}^g_{dd}  :=  \left\{%
        \begin{array}{l}  
            f:C \rightarrow {\mathbb P}^1 \ \mbox{such that}\\
            \ \ \ \bullet\ C\ \mbox{is a connected  curve of genus $g$;}\\
          \ \ \   \bullet f \ \mbox{is a map of degree $d$;} \\
          \ \ \   \bullet f \ \mbox{is fully ramified over $0$ and $\infty$;} \\
          \ \ \  \bullet f \ \mbox{is (at most) simply ramified over all other points} \\
        \end{array} \right\} /PGL(2,\mathbb C)
$$
and their induced source maps $s:{\mathcal{A}}_{dd}^g \rightarrow {\mathcal{M}_g}$ to the moduli
space of curves(there is some sloppiness in this definition, as one needs to also factor by the automorphisms
of $f$).
The Hurwitz spaces determine classes $[{\mathcal{A}}_{dd}^g] \in R^{g-2}({\mathcal{M}_g})$ 
that are multiples of a single class (\cite{l:ottr}).

\medskip

There is a smooth compaticfication $\mathcal{A}^g_{dd} \subset \overline {\mathcal{A}}_{dd}^g$ by the 
stack of admissible 
covers (section \ref{la})  such that $s$ extends to a map 
$ s:\overline {\mathcal{A}}_{dd}^g \rightarrow \overline {\mathcal{M}_g}$ to the Deligne-Mumford 
moduli stack of stable curves. A Hodge bundle  $\mathbb{E}$ lives on $\overline {\mathcal{M}_g}$,
whose Chern classes $\lambda_1,\dots \lambda_g$ are all tautological.
Thus the following generating functions (if non-zero!) completely determine the proportionalities  of the 
tautological classes
$[\mathcal{A}_{dd}^g] \in R^*({\mathcal{M}_g})$:
$${\mathcal D}(u) := \sum_{g = 1}^\infty \left(\int_{\overline {\mathcal{A}}_{dd}^g} \lambda_g \lambda_{g-1}
\right) \frac{u^{2g}}{2g!}.$$

\textbf{Remark:} It is not difficult to show that the class $\la$ vanishes on $\overline{A}_{dd}^g \smallsetminus {\mathcal{A}}_{dd}^g$. This allows one to regard capping with $\la$ as an evaluation on the tautological ring $R^\ast(\mathcal{M}_g)$.

These generating functions were first computed by Bryan-Pandharipande in \cite{bp:tlgwtoc}, 
using Hodge integral techniques,
a sophisticated geometric argument from Looijenga, and Faber's use of Mumford's computations in 
the hyperelliptic case $(d=2)$ that made use of the Grothendieck-Riemann-Roch theorem(!):

\begin{theorem}[\cite{bp:tlgwtoc}, Thm 6.5] 
\label{bpa}
$${\mathcal D}(u) = \ln\left(\frac{d\sin\left(\frac{u}{2}\right)}{\sin\left(\frac{du}{2}\right)}\right)$$
\end{theorem}
\medskip

Our point in this paper is to show how this (and other related) generating functions follow directly from
localization on spaces of admissible covers. This reduces the computation to the combinatorics of simple 
Hurwitz numbers and, perhaps equally importantly, reduces the proof to a simple tallying up of fixed loci 
and their contributions to an equivariant integral. In other words, through localization, all the geometry has essentially
been reduced to combinatorics.
 
\medskip
 
\noindent {\bf Acknowledgement:} Some suspicious-looking 
similarities among generating functions 
in local Gromov-Witten theory prompted
Ravi Vakil to demand an explanation of the second author. This paper (we hope)
serves as an explanation.

\medskip

\section{The Generating Functions}
This section serves the dual purpose of establishing our notation and generating function conventions, and introducing the main characters in the forthcoming computations. When dealing with families of enumerative invariants, it is very important to organize the information both naturally and
efficiently. Our generating function conventions are driven both by geometry 
 (in the choice of exponents of the 
formal parameters) and 
combinatorics  (in the simple closed formulas for the generating functions with these conventions).
This convention is already on display in the generating function ${\mathcal D}(u)$ of the 
introduction, where the 
exponent $u^{2g}$ and denominator $2g!$ correspond to the number of simple ramification points
($2g$) of each of the maps $f:C \rightarrow \mathbb P^1$ of $\mathcal{A}^g_{dd}$.

	\subsection{Simple Hurwitz Numbers}
	\label{hurwitz}
	Simple Hurwitz numbers count ramified covers of the Riemann sphere with prescribed ramification data over one special point of $\proj$ and simple ramification over the other (fixed) branch points. 
	Let $\eta$ be a partition of the integer $d$. Then the
associated simple Hurwitz number in genus $g$ is defined to be:
$$H^g_\eta := \sum_{\pi \in S^g_\eta} \frac 1{|\mbox{Aut}(\pi)|}$$
where $S^g_\eta$ is the following set:
$$
    \begin{array}{ccc}
        S^g_\eta & := &\ \left\{%
        \begin{array}{l}
            \mbox{degree  $d$  covers} \\
            \  \ \ \ \ \ \ \ \ \ \ \ \ \ \ \ \ \ {C \stackrel{\pi}{\longrightarrow} \proj}\\ 
            \mbox{such that:}  \\
            \bullet\ C\ \mbox{is a connected  curve of genus $g$;}\\
            \bullet \ \{p_1,\ldots, p_r\} \mbox{\ are marked points in}\\ \mbox{fixed generic position on $\proj$;}\\
            \bullet\ \pi \ \mbox{is  unramified over}  \\
            \ \ \ \ \ \ \ \ \ \ \ \ \ \proj \smallsetminus\{\infty,p_1,\ldots, p_r\};\\
            \bullet \ \pi\ \mbox{ramifies with profile $\eta$ over $\infty$;} \\ 
            \bullet\ \pi\ \mbox{has simple ramification over $p_i$.} 
        \end{array}
        \right\}
    \end{array}
$$
and Aut$(\pi)$ is the (finite!) automorphism group of each such map $\pi$.

\medskip

\noindent 

\textbf{Remark:} The number $r$ of simple branch points is determined by the Riemann-Hurwitz formula:
$$r=2g+d-2+\ell(\eta).$$
where $\ell(\eta)$ is the length of the partition $\eta$.

\medskip

Using this as a definition for $r$, we organize Hurwitz numbers in generating function form as follows:

\begin{eqnarray}
\hf{\eta}(u):= \sum_{g=0}^\infty (-1)^g\hn{g}{\eta}\ \frac{u^r}{r!}
\end{eqnarray}

\noindent 
\textbf{Remark:} The sign $(-1)^g$ is non-standard, and has no particular meaning but to tune this generating function to our localization set-up.

\medskip

For the $d$-cycle $\eta = (d)$, the generating function $\hf{(d)}(u)$ has a particularly simple form, that was first computed by Shapiro, Shapiro and Vainshtein in \cite{ssv:rc}:

\begin{formula}[\cite{ssv:rc}, Thm 6]
\label{Hd}
\begin{eqnarray*}
\begin{imp}
	\displaystyle{\hf{(d)}(u)= \frac{2^{d-1}}{d\cdot d!}\left(\sin\left(\frac{du}{2}\right)\right)^{d-1}.}
\end{imp}
\end{eqnarray*}
\end{formula} 

In general, the single Hurwitz numbers are harder to compute. However, the following relation 
among single Hurwitz numbers
was obtained by localization techniques by the second author in \cite{r:tqft}:

\begin{formula}[\cite{r:tqft}, Equation (26)]
\label{Heta}
\begin{eqnarray*}
\begin{imp}
	\displaystyle{0 = \sum_{\eta \vdash d} (-1)^{\ell(\eta)} \frac
	{\hf{\eta}(u)}{\prod_{\eta_i\in \eta} 2\sin\left( \frac{\eta_iu}2\right)}}
\end{imp}
\end{eqnarray*}
where $\eta_i \in \eta$ are the cycle-lengths of the partition $\eta$.
\end{formula} 

\subsection{$\lambda_g\lambda_{g-1}$ Hurwitz-Hodge Integrals}
\label{la}
Let:
$$\overline {\mathcal{A}}^g_{dd} =\A{g}{0}$$
 be the Abramovich-Corti-Vistoli (smooth, proper) moduli stack (\cite{acv:ac}) of:
\begin{itemize}
	\item genus $g$, degree $d$ admissible covers of a genus $0$ (possibly nodal) curve;
	\item with all $2g+2$ branch points marked;
	\item the first $2g$ branch points carrying a simple ramification profile;
	\item the last two branch points carrying ramification profile $(d)$.
\end{itemize}

The universal family of admissible covers defines source and target maps:
$$
\begin{array}{ccc}
 \overline{\mathcal{A}}^g_{dd} & \stackrel{s} \rightarrow & \overline{\mathcal{M}_g} \\
\ \ \ \downarrow t \\ \overline{\mathcal{M}}_{0,2g+2} 
\end{array}
$$
to the moduli spaces of  genus $g$ and genus $0$, $2g+2$-pointed stable curves.

\medskip

Let $\lambda_i$ be the Chern classes of the Hodge bundle (pulled back from $\overline {\mathcal{M}_g}$), 
and in addition, let 
$\psi = t^*\psi_{2g+2}$ be the pull-back of the Mumford-Morita class from $\overline{\mathcal{M}}_{0,2g+2}$ (\cite{m:taegotmsoc}). 
Integrals of products of such classes are referred to as Hurwitz-Hodge integrals (\cite{bgp:crc}) and 
following \cite{r:dd2}, we define the following generating functions of Hurwitz-Hodge integrals:
\begin{eqnarray}
\mathcal{D}_i(u):=\sum_{g\geq i} \left( \int_{\overline {A}^g_{dd}} \lambda_g\lambda_{g-i}\psi^{i-1} \right) \frac{u^{2g}}{2g!}	
\end{eqnarray}
In particular, ${\mathcal D}(u) = {\mathcal D}_1(u)$.

\medskip

These generating functions are all related via:
\begin{formula}[\cite{r:dd2}, Thm 1]
\label{dlf}
\begin{eqnarray*}
\begin{imp}
	\displaystyle{\mathcal{D}_i(u):= \frac{d^{i-1}}{i!}\left(\mathcal{D}(u)\right)^i}
\end{imp}
\end{eqnarray*}
\end{formula}

\noindent 
\textbf{Remark:} Thus if we set ${\mathcal D}_0(u) := \frac 1d$, then we may combine the
${\mathcal D}_i(u)$ in:

\begin{eqnarray}
\label{cor}
	\displaystyle{\mathcal{T}(u)=\sum_{g\geq 0} T^g \frac{u^{2g}}{2g!}:=\sum_{i=0}^\infty \mathcal{D}_i(u)= 
	\frac{1}{d}\ \mbox{e}^{d\mathcal{D}(u)}.}
\end{eqnarray}

\subsection{CY Cap}
\label{CYcap}
Consider the integral:

\begin{eqnarray}
	 CY^g:=\int_{\star} c_{g}(\mathbb{E}^\ast)c_{g+d-1}( R^1\pi_\ast f^\ast(\mathcal{O}_{\proj}(-1)))ev^\ast_{(d)}(\infty),
\label{int}
\end{eqnarray}

where:

\begin{itemize}
	\item 
$$
\star:=\overline{Adm}_{g\stackrel{d}{\rightarrow}\proj,((d),t_1,\ldots,t_{2g+d-1})}
$$
is the moduli space of degree $d$, genus $g$ admissible covers of a parametrized $
\proj$, fully ramified over one special branch point, simple ramification over all other (marked!) branch points. (see \cite{r:thesis}, section 1.1.6.2)
	\item $\pi: \mathcal{U} \rightarrow \star$ is the projection from the universal admissible cover family;
	\item $f:\mathcal{U}\rightarrow \proj$ is the universal map \footnote{followed by contraction of the ``sprouted twigs''.};
	\item $ev_{(d)}: \star \rightarrow \proj$ is the map evaluating a cover at its unique fully ramified point (i.e. the composition of $f$ with the first section $s_{(d)}: \star \rightarrow \mathcal{U}$).   
\end{itemize}

Define the generating function:
\begin{eqnarray}
\mathcal{CY}(u) := \sum_{g=0}^\infty CY^g\ \frac{u^{2g+d-1}}{({2g+d-1})!}
\end{eqnarray}	 
Localization and relations among simple Hurwitz numbers determine this generating function.
The computation is carried out in \cite{r:adm} for degrees $2,3$ and in \cite{r:tqft} for  general degree. The following result is
obtained by localizing and reducing to Formula \ref{Heta} above:
\begin{formula}[\cite{r:tqft}, Thm 3]
\label{CY}
\begin{eqnarray*}
\begin{imp}
	\displaystyle{\mathcal{CY}(u)= {(-1)^{d-1}}\frac{1}{d}\frac{\left(2\sin\left(\frac{u}{2}\right)\right)^d}{2\sin\left(\frac{du}{2}\right)}.}
\end{imp}
\end{eqnarray*}
\end{formula}
\textbf{Remarks:}
\begin{enumerate}
	\item The generating function $CY$ is a fundamental quantity in the TQFT of Admissible Cover invariants (\cite{r:tqft}): it represents the only non-zero connected invariants for the Calabi-Yau cap, i.e. the one pointed invariants for the zero section of the bundle $\mathcal{O}_{\proj} \oplus \mathcal{O}_{\proj}(-1)$. It is also a close relative of the analogous invariants in the local GW theory of curves (\cite{bp:tlgwtoc}). This explains its otherwise slightly mysterious presence (and name) in this application.
	\item There is one small notational difference between the description of $\mathcal{CY}(u)$ here and in \cite{r:tqft}. There, we worked on moduli spaces of covers where the branch locus on $\proj$ was not marked. Here we  mark all branch points on $\proj$, and then include the ``normalization factor" $1/({2g+d-1})!$ to compensate. The two choices are obviously equivalent.
\end{enumerate}
\section{Proof of theorems}
	\subsection{Localization Computation}
\subsubsection{The Set-up}

Consider the  one-dimensional algebraic torus $\mathbb{C}^\ast$, and recall that the $\mathbb{C}^\ast$-equivariant Chow ring of a point is a polynomial ring in one variable:
$$A^\ast_{\Cstar}(\{pt\},\mathbb{C})= \mathbb{C}[\hbar]. $$
Let $\mathbb{C}^\ast$ act on a two-dimensional vector space $V$ via:
$$t\cdot(z_0,z_1)=(tz_0,z_1).$$
This action descends to $\proj$, with fixed points $0=(1:0)$ and $\infty=(0:1)$. An equivariant lifting of $\Cstar$ to a line bundle $L$ over $\proj$ is uniquely determined by its  weights $\{L_0,L_\infty\}$ over the fixed points.

The action on $\proj$ induces an action on the moduli spaces of admissible covers to a parametrized $\proj$ simply by post composing the cover map with the automorphism of $\proj$ defined by $t$. 

The fixed loci for the induced action on the moduli space consist of admissible covers such that anything ``interesting'' (ramification, nodes, marked points) happens over $0$ or $\infty$, or on ``non special'' twigs that attach to the main $\proj$ at  $0$ or $\infty$.	

\subsubsection{The Strategy}
\label{strategy}
	In \cite{r:tqft}, the second author proved Formula \ref{CY} using localization. The integrand in (\ref{int}) can be regarded as an equivariant cohomology class. A well chosen linearization for the bundles $\mathcal{O}_{\proj}$ and $\mathcal{O}_{\proj}(-1)$ allows the explicit evaluation of it - which is a non-equivariant quantity. Hence, such evaluation is independent of the linearization chosen.  Here we use a different system of linearization weights to evaluate (\ref{int}):
	\begin{center}
\begin{tabular}{|l||c|c|}
\hline
      & weight over $0$ & weight over $\infty$ \\
\hline
\hline
 ${\mathcal{O}_{\proj}}(-1)$ & -1 & 0 \\
\hline
${\mathcal{O}_{\proj}}$  & 1  & 1  \\
\hline
\end{tabular}
\end{center}
\subsubsection{The possibly contributing fixed loci}
\begin{figure}[htbp]
	\begin{center}
$$	
\begin{array}{ccc}
	F_g = \ \ 	\includegraphics[width=0.32\textwidth]{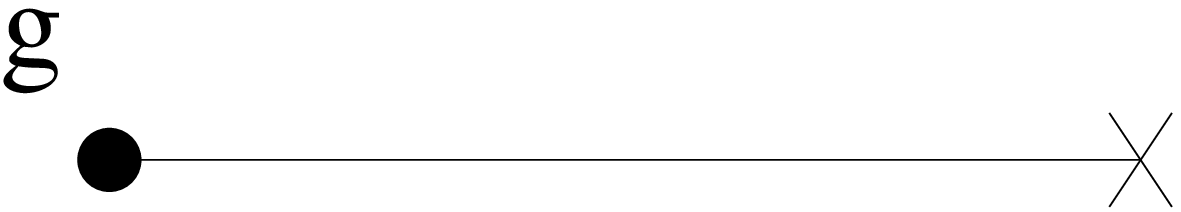} &
				\hspace{.7cm}																							&		
	F_{g_1g_2}=\  \	\includegraphics[width=0.35\textwidth]{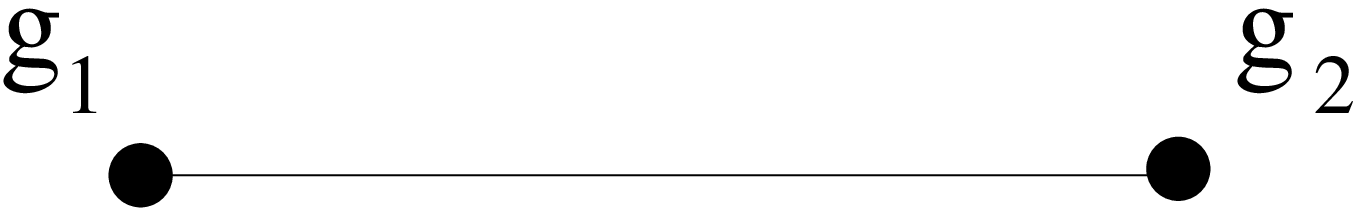}
\end{array}
$$
	\end{center}
	\caption{The localization graphs of the possibly contributing fixed loci.}
	\label{fix}
\end{figure}

The above choice of linearization restricts the number of possibly contributing fixed loci (identified here with their localization graphs):
\begin{enumerate} 
	\item The full ramification condition  at $\infty$ forces the preimages of  $\infty$ to be connected. 
	\item The weight ``$0$" over $\infty$ forces the localization graphs to have valence $1$ over $\infty$ (see, e.g.,  \cite{r:thesis}).
\end{enumerate}
The possibly contributing fixed loci consist of boundary strata parameterizing a single sphere, fully ramified over $0$ and$\infty$, mapping with degree $d$ to the main $\proj$. A rational tail  $T_0$ must sprout from $0$, and be covered by a curve of genus $g_1$, $0\leq g_1 \leq g$. If $g_1<g$, then a rational tail $T_\infty$ sprouts from $\infty$ as well, covered by a curve $C$ of genus $g_2=g-g_1$. The map $C \rightarrow T_\infty$ has two points of full ramification: the attaching point and another point, corresponding to the original mark $(d)$. 
The situation is illustrated in Figure \ref{fix}.

These fixed loci are covers  of products of spaces of admissible covers, with ``stacky" multiplicities (see \cite{r:thesis} for a quick explanation, or \cite{agv:gwtodms} for a comprehensive treatise of this subject):
$$
\begin{array}{lcl}
	 \displaystyle{F_g\cong  \overline{\mathcal{A}_{d}^g}} & \hspace{2cm} & 
	 \displaystyle{F_{g_1g_2}\rightarrow   \overline{\mathcal{A}}_{d}^{g_1} \times \overline{\mathcal{A}}_{dd}^{g_2}} \mbox{\ of degree $d {{2g+d-1}\choose{2g_2}}$},
\end{array}
$$
with the following notation:

\begin{itemize}
	\item $\overline{\mathcal{A}}^g_d := \As{g}{0}$;
	\item $\overline{\mathcal{A}}^g_{dd} := \A{g}{0}$.
\end{itemize}

	\subsubsection{The relations}
Before we write down the explicit contributions to integral (\ref{int}), let us introcuce one last piece of notation that is intended to make the computation  a little less painful to read:

 $$\Lambda_g(n):=\sum_{i=0}^g (n\hbar)^i\lambda_{g-i}.$$
	
With the choice of linearization from Section  \ref{strategy}, the contribution to $CY^g$ from the fixed locus $F_{g_1g_2}$ is:
\begin{eqnarray*}
	d\binom {2g+d-1}{2g_2}(-\hbar)^{d-1}\frac{(d-1)!}{d^{d-1}}\int_{{\overline{\mathcal{A}}}^{g_1}_d}\frac{\Lambda_{g_1}(-1)\Lambda_{g_1}(+1)}{\hbar(\hbar-\psi)}\int_{{\overline{\mathcal{A}}}^{g_2}_{dd}}\frac{\Lambda_{g_2}(-1)\lambda_{g_2}}{\hbar(\hbar+\psi)}. \nonumber
\end{eqnarray*}
	Expanding the denominators and recalling Mumford's relation $c(\mathbb{E}\oplus \mathbb{E}^\ast)=1$, this reduces to:
\begin{eqnarray}
\label{fgg}
\frac{(-1)^{d-1}d!}{d^{d-1}}\binom {2g+d-1}{2g_2}\int_{{\overline{\mathcal{A}}}^{g_1}_d}(-1)^{g_1}\psi^{2g_1+d-3} 
\int_{{\overline{\mathcal{A}}}^{g_2}_{dd}}\lambda_{g_2}\lambda_{g_2-1}+\ldots+\lambda_{g_2}\psi^{g_2-1}= \nonumber\\
=\frac{(-1)^{d-1}d!}{d^{d-1}}(2g+d-1)!\ H_{(d)}^{g_1}T^{g_2}. 
\end{eqnarray}	
	The contribution of $F_g$ is:
	
\begin{eqnarray*}
	(-\hbar)^{d-1}\frac{(d-1)!}{d^{d-1}}\int_{{\overline{\mathcal{A}}}^{g}_d}\frac{\Lambda_{g}(-1)\Lambda_{g}(+1)}{\hbar(\hbar-\psi)}
	= \frac{(-1)^{d-1}(d-1)!}{d^{d-1}}(2g+d-1)!\ H_{(d)}^{g},
\end{eqnarray*}
which, since $T^0=\mathcal{D}_0=1/d$, can be rewritten as
\begin{eqnarray}
\label{Fg}
	\frac{(-1)^{d-1}d!}{d^{d-1}}(2g+d-1)!\ H_{(d)}^{g}T^{0}.
\end{eqnarray}

Equations (\ref{fgg}) and (\ref{Fg}) are nicely encoded in generating function form to give the relation:
\begin{eqnarray}
	\label{relation}
	\mathcal{CY}(u)=\frac{(-1)^{d-1}d!}{d^{d-1}} \hf{(d)}(u) \mathcal{T}(u).
\end{eqnarray}

Theorem \ref{rel} follows immediately from (\ref{relation}) and  (\ref{cor}).

Theorem \ref{bpa} is a straightforward corollary of  Theorem \ref{rel}:
recall that the combinatorics of simple Hurwitz numbers gives us explicit expressions for $\hf{(d)}(u)$ (Formula   \ref{Hd}) and $\mathcal{CY}(u)$ (Formula \ref{CY}). Solving for $\mathcal{D}(u)$ we obtain:


$$
\begin{imp}
\label{fapalo}
\displaystyle{
	\mathcal{D}(u)=\ln\left(\frac{d\sin\left(\frac{u}{2}\right)}{\sin\left(\frac{du}{2}\right)}\right).
	}
\end{imp}
$$

\addcontentsline{toc}{section}{Bibliography}
\bibliographystyle{alpha}
\bibliography{biblio}

\end{document}